\documentclass[10pt,twoside,a4paper,reqno]{amsart}
\usepackage{amsfonts,amsthm,latexsym,amsmath,amssymb,amscd,epsf}

\theoremstyle{plain}
\newtheorem{theo}           {Theorem}
\newtheorem{pro}            {Proposition}
\newtheorem{coro}           {Corollary}
\newtheorem{lemm}           {Lemma}
\newtheorem{conj}           {Conjecture}

\theoremstyle{definition}
\newtheorem{pr}              {Problem}
\newtheorem*{ack}            {Acknowledgements}

\theoremstyle{remark}
\newtheorem{rem}             {Remark}

\newenvironment{theorem}{\begin{theo}}{\end{theo}}

\newenvironment{corollary}{\begin{coro}}{\end{coro}}
\newenvironment{lemma}{\begin{lemm}}{\end{lemm}}
\newenvironment{remark}{\begin{rem}}{\end{rem}}

\newenvironment{problem}{\begin{pr}}{\end{pr}}

\newcommand \Ga {\Gamma}
\newcommand \ga {\gamma}

\newcommand{\om}{\omega}
\newcommand{\ze}{\zeta}

\newcommand \bR {\mathbb R}

\newcommand \bZ {\mathbb Z}
\newcommand\al {\alpha}
\newcommand \be {\beta}

\begin{document}

\title[Root asymptotics of spectral polynomials]
{Root asymptotics of spectral polynomials \\ for the  Lam\'e operator}

\author[J.~Borcea]{Julius Borcea}
\address{Department of Mathematics, Stockholm University, SE-106 91
Stockholm,
      Sweden}
\email{julius@math.su.se}

\author[B.~Shapiro]{Boris Shapiro}
\address{Department of Mathematics, Stockholm University, SE-106 91
Stockholm,
      Sweden}
\email{shapiro@math.su.se}

\begin{abstract}
  The study of polynomial solutions to the classical Lam\'e equation in its
algebraic form, or equivalently,
  of double-periodic solutions of its Weierstrass form has a long history.
Such solutions appear at integer values of the spectral parameter and
their respective eigenvalues serve as the ends of bands in the boundary value
problem for the corresponding Schr\"odinger equation with  finite gap
potential given by the Weierstrass $\wp$-function on the real line.
In this paper we establish several natural (and equivalent) formulas in terms
of hypergeometric and elliptic type integrals for the density of the
appropriately scaled asymptotic distribution of these eigenvalues when the
integer-valued spectral parameter tends to infinity. We also show that this
density satisfies a Heun differential equation with four singularities.
\end{abstract}

\maketitle

\section{Introduction and Main Results}

The algebraic form of the classical Lam\'e equation \cite[Chap.~23]{WW} was
introduced by Lam\'e in the 1830's in connection with the  separation of
variables in the Laplace equation by means of elliptic coordinates in $\bR^l$.
Lam\'e's equation is given by
\begin{equation}\label{eq:standLame}
\left\{Q(z)\frac{d^2}{dz^2}+\frac{1}{2}Q'(z)\frac{d}{dz}+V(z)\right\}S(z)=0,
\end{equation}
where $Q_l(z)$ is a real degree $l$ polynomial with all real and distinct
roots and $V(z)$ is a polynomial of degree at most $l-2$ whose choice depends
on the type of solution to \eqref{eq:standLame} one is looking for.
In the second half of the nineteenth century several famous mathematicians
including B\^ocher, Heine, Klein and Stieltjes studied the number and
various properties of the so-called {\em Lam\'e solutions of the first kind}
(of given degree and type) to equation \eqref{eq:standLame}.
These are also known as {\em Lam\'e polynomials} of a certain type.
Such solutions exist for special choices of $V(z)$
and are characterized by the property that their logarithmic derivative is a
rational function. For a given $Q(z)$ of degree  $l\ge 2$ with simple roots
there exist $2^l$ different possibilities for Lam\'e solutions depending on
whether these solutions are smooth at a given root of $Q(z)$ or have  a
square root singularity, see \cite{Poo,WW} for more details.

A {\em generalized Lam\'e equation} \cite{WW} is a second
order differential equation of the form
\begin{equation}\label{eq:comLame}
         \left\{Q(z)\frac
         {d^2}{dz^2}+P(z)\frac{d}{dz}+V(z)\right\}S(z)=0,
         \end{equation}
where $Q(z)$ is a complex polynomial of degree $l$ and $P(z)$ is
a complex polynomial of degree at most $l-1$.
As it was first shown by Heine \cite{He} for a generic equation of either of
the forms \eqref{eq:standLame} or
\eqref{eq:comLame} and an arbitrary positive integer $n$ there are exactly
$\binom{n+l-2}{n}$ polynomials $V(z)$ such that $S(z)$ is a polynomial of
degree $n$. Below we concentrate on the most classical case of a
cubic polynomial $Q(z)$ with all real roots as treated in e.g.~\cite{WW}.
Already Lam\'e and Liouville knew that in this case the (unique) root of
$V(z)$ is real and located between the minimal and the maximal roots of
$Q(z)$. The latter property was further generalized by Heine, Van Vleck and
other authors. Note that if $\deg  Q(z)=3$ then $V(z)$ is at most linear and
that for a given value of the positive integer $n$ there are at most $n+1$
such polynomials.

Before formulating our main results let us briefly review -- following
mostly \cite[\S 3]{Ma2} and \cite{Poo} -- some necessary background on
the version of the Lam\'e equation used below and its spectral polynomials.
Setting, as one traditionally does, $Q(z)=4(z-e_1)(z-e_2)(z-e_3)$ with
$e_1> e_2> e_3$ we can rewrite equation \eqref{eq:standLame} as
\begin{equation}\label{eq:myLame}
\left\{ \frac{d^2}{dz^2}+\frac{1}{2}\sum_{i=1}^3\frac{1}{z-e_i}\frac{d}{dz}
-\frac{n(n+1)z+E}{4\prod_{i=1}^3(z-e_i)}\right\}S(z)=0.
\end{equation}
(The chosen representation for the linear polynomial $V(z)$ will shortly become
clear.) Notice that several equivalent forms of the Lam\'e equation are
classically known. Among those one should mention two algebraic forms, the
Jacobian form and the Weierstrassian form, respectively,
see \cite[\S 23.4]{WW} and Remark \ref{r-wei} below.
Equation \eqref{eq:myLame} presents the most commonly used real algebraic
form of the Lam\'e equation which is smooth while the other (real) algebraic
form is singular.\footnote{The authors thank the anonymous referee for this
observation.}

A {\em Lam\'e solution of the first kind} to equation
\eqref{eq:myLame} is a solution of the form
$S(z)=(z-e_1)^{\kappa_1}(z-e_2)^{\kappa_2}(z-e_3)^{\kappa_3}\tilde S(z)$,
where each $\kappa_i$ is either $0$ or $\frac{1}{2}$ and $\tilde S(z)$ is a
polynomial. Lam\'e solutions which are pure polynomials -- i.e., for
which $\kappa_1=\kappa_2=\kappa_3=0$ -- are said to be of
{\em type 1},  those with
a single square root are of {\em type 2}, those with two square roots are of
{\em type 3}  and, finally,
those involving three square roots are said to be of {\em type 4}.
(``Types'' are sometimes called ``species'', see, e.g., \cite {WW}.)
One can easily check that \eqref{eq:myLame} has a Lam\'e solution if and only
if  $n$ is a nonnegative integer. Moreover, if $n$ is even then only
solutions of types $1$ and $3$ exist. Namely, for appropriate choices of the
energy constant $E$ one gets exactly $\frac{n+2}{2}$ distinct independent
solutions of type $1$ and $\frac{3n}{2}$ independent solutions of type $3$.
If $n$ is odd then only solutions of types $2$ and $4$ exist. In this case,
for appropriate choices of the energy constant $E$ one gets exactly
$\frac{3(n+1)}{2}$ distinct independent solutions of type $2$ and
$\frac{n-1}{2}$ independent solutions of type $4$. In both cases the
total number of Lam\'e solutions equals $2n+1$, which coincides with the
number of independent spherical harmonics of order $n$.

Let $R_n(E)=\prod_{j=0}^{2n+1}(E-E_j)$, $n\in\bZ_{+}$, denote the monic
polynomial of degree $2n+1$ whose roots are exactly the values of the energy
$E$ at which equation \eqref{eq:myLame} has a Lam\'e solution. These
polynomials are often referred  to as {\em spectral polynomials} in the
literature; their
study goes back to Hermite  and Halphen. The most recent results in this
direction can be found in \cite{BeEn,Ma2,Ta1,Ta2,Vo} and \cite{GrVe1,GrVe2}.
In particular, \cite{Ma2} contains an excellent survey of this topic as well
as a comprehensive table with these polynomials (and their modified versions)
correcting several mistakes that occurred in previous publications. Article
  \cite{GrVe2} is apparently the first attempt to give a (somewhat) closed
formula for $R_n(E)$; for this the authors use yet another
family of polynomials which they call {\em elliptic Bernoulli polynomials}.
Since explicit formulas for $R_n(E)$ seem to be rather complicated  and
difficult to handle, in this paper we  study the asymptotics of the root
distribution of appropriately scaled versions of $R_n(E)$
(Corollary~\ref{cor:main}). In spite of a more than
150 years long history of the Lam\'e equation the only source discussing
questions similar to ours that we were able to locate is \cite{Bo}. We should
also mention that the results below are actually much more precise and
therefore supersede those of {\em loc.~cit.}

Let us now introduce eight spectral  polynomials
$R^{\kappa_1,\kappa_2,\kappa_3}_n(E)$ related to the eight  types of
Lam\'e solutions mentioned above, namely
\begin{equation} \label{eq:spect}
  R_n(E)=\begin{cases}R^{0,0,0}_n(E)R^{\frac{1}{2},\frac{1}{2},0}_n(E)
R^{\frac{1}{2},0,\frac{1}{2}}_n(E)R^{0,\frac{1}{2},\frac{1}{2}}_n(E)
\text{ when } n \text{ is even,}\\
  R^{\frac{1}{2},0,0}_n(E)R^{0,\frac{1}{2},0}_n(E)R^{0,0,\frac{1}{2}}_n(E)
R^{\frac{1}{2},\frac{1}{2},\frac{1}{2}}_n(E)  \text{ when } n \text{ is odd.}
  \end{cases}
  \end{equation}
Note for example that if $n$ is even then $R^{0,0,0}_n(E)$ is the (unique)
monic polynomial
of degree $\frac{n+2}{2}$ whose roots are precisely the values of $E$
for which equation \eqref{eq:myLame} has a pure polynomial solution.

The results of the present paper actually hold not only for equation
\eqref{eq:myLame} but also for generalized Lam\'e equations
(cf.~\eqref{eq:comLame}) of the form
\begin{equation}\label{eq:myGenLame}
\left\{ \frac{d^2}{dz^2}
+\sum_{i=1}^3\frac{\al_i}{z-e_i}\frac{d}{dz}
+\frac{V(z)}{\prod_{i=1}^3(z-e_i)}\right\} S(z)=0,
\end{equation}
where $\al_i>0$, $i=1,2,3$, and $V(z)$ is an undetermined linear polynomial.
This case was thoroughly treated by Stieltjes in \cite{St}. In particular,
he proved that for any positive integer $n$ there exist exactly $n+1$
polynomials $V_{n,j}(z)$, $1\le j\le n+1$, such that \eqref{eq:myGenLame}
has a polynomial solution $S(z)$ of degree $n$. Moreover, the unique
root $t_{n,j}$ of $V_{n,j}(z)$ lies in the interval $(e_3,e_1)$ and these
$n+1$ roots are pairwise distinct. Consider now the polynomial
$$Sp_n(t)=\prod_{j=1}^{n+1}(t-t_{n,j}),$$
which is the scaled version of the spectral polynomial
$R_n^{0,0,0}(E)$. More exactly, for any even $n$ one has that 
$Sp_n(t)=\frac{R_n^{0,0,0}(n(n+1)t)}{(n(n+1))^\frac{n+2}{2}}$. Assume 
further that the set $\{t_{n,j}:1\le j\le n+1\}$
is ordered so that $t_{n,1}<t_{n,2}<\ldots<t_{n,n+1}$ and associate to it
the finite measure
$$\mu_n=\frac{1}{n+1}\sum_{j=1}^{n+1}{\delta(z-t_{n,j})},$$
where $\delta(z-a)$ is the Dirac measure supported at $a$.
The measure $\mu_n$ thus obtained is clearly a real probability measure that
one usually refers to as the root-counting measure of $Sp_n(t)$.

Below we shall make use of some well-known notions from
the theory of special functions that can be found in e.g.~\cite{AS}. In
fact several of our formulas rely on various integral
representations and transformation properties of the Gauss hypergeometric
series
$$F(a,b,c;z)=_{2}\!\!F_{1}(a,b,c;z)=\frac{\Ga(c)}{\Ga(a)\Ga(b)}
\sum_{m=0}^{\infty}\frac{\Ga(a+m)\Ga(b+m)}{\Ga(c+m)}\frac{z^m}{m!},
$$
where $\Ga$ denotes Euler's Gamma-function, as well as the complete
elliptic integral of the first kind
\begin{equation}\label{eq:elliptic}
\mathbb K(\zeta)=\int_0^1\frac{dt}{\sqrt{(1-t^2)(1-\zeta^2t^2)}}.
\end{equation}
Recall from \cite[\S 15.3.1]{AS} that if $\Re(c)>\Re(b)>0$ then a convenient
way of rewriting $F(a,b,c;z)$ is given by
\begin{equation}\label{eq:hyper}
F(a,b,c;z)=\frac{\Ga(c)}{\Ga(b)\Ga(c-b)}
\int_{0}^{1}t^{b-1}(1-t)^{c-b-1}(1-tz)^{-a}dt.
\end{equation}
The above integral represents a one-valued analytic function in the $z$-plane
cut along the real axis from $1$ to $\infty$ and thus it gives the analytic
continuation of $F(a,b,c;z)$. To simplify our formulas let us also define a
real-valued function of three real variables
$(x_1,x_2,x_3)\mapsto \om(x_1,x_2,x_3)$ by setting
\begin{equation}\label{eq:om}
\om(x_1,x_2,x_3)=\frac{(|x_1|+|x_2|)|x_3|}{|x_1||x_2-x_3|+|x_2||x_1-x_3|}.
\end{equation}
Note that $\om(x_1,x_2,x_3)$ is symmetric in its first two arguments and that
\begin{equation}\label{eq:omega}
0\le \om(x_1,x_2,x_3)\le 1\text{ whenever }x_1<0,\,x_2>0\text{ and }
x_1\le x_3\le x_2.
\end{equation}
We should emphasize that the latter conditions on $x_1,x_2,x_3$ -- hence also
the bounds for $\om(x_1,x_2,x_3)$ given in \eqref{eq:omega} -- always hold in
the present context.

We are now ready to state our main results.

\begin{theorem}\label{th:main}
For any choice of real numbers $e_1>e_2>e_3$ and positive numbers
$\al_1,\al_2,\al_3$ the sequence of measures
$\{\mu_n\}_{n\in\bZ_{+}}$ strongly converges to the probability measure
$\mu_Q$ supported on the interval $[e_3,e_1]$ with density $\rho_Q$ given by
any of the following equivalent expressions:
\begin{itemize}
\item[(i)]
\begin{multline*}
\rho_Q(s)=\frac{1}{\pi}
\sqrt{\frac{1+\sqrt{1-\om(e_1-e_2,e_3-e_2,s-e_2)^2}}
{(e_1-e_3)|s-e_2|\om(e_1-e_2,e_3-e_2,s-e_2)}}\\
\times\mathbb K\!\left(i\sqrt{\frac{2\sqrt{1-\om(e_1-e_2,e_3-e_2,s-e_2)^2}}
{1-\sqrt{1-\om(e_1-e_2,e_3-e_2,s-e_2)^2}}}\right),
\end{multline*}
where $\om$ is the function defined in \eqref{eq:om}.
\item[(ii)]
\begin{multline*}
\rho_Q(s)=\frac{1}{2}\sqrt{\frac{1+\sqrt{1-\om(e_1-e_2,e_3-e_2,s-e_2)^2}}
{(e_1-e_3)|s-e_2|\om(e_1-e_2,e_3-e_2,s-e_2)}}\\
\times F\!\left(\frac{1}{2},\frac{1}{2},1;
-\frac{2\sqrt{1-\om(e_1-e_2,e_3-e_2,s-e_2)^2}}
{1-\sqrt{1-\om(e_1-e_2,e_3-e_2,s-e_2)^2}}\right).
\end{multline*}
\item[(iii)]
\begin{multline*}
\rho_Q(s)=\sqrt{\frac{\om(e_1-e_2,e_3-e_2,s-e_2)}{2(e_1-e_3)|s-e_2|
(1+\om(e_1-e_2,e_3-e_2,s-e_2))}}\\
\times F\!\left(\frac{1}{2},\frac{1}{2},1;
\frac{1-\om(e_1-e_2,e_3-e_2,s-e_2)}{1+\om(e_1-e_2,e_3-e_2,s-e_2)}\right).
\end{multline*}
\item[(iv)]
\begin{multline*}
\rho_Q(s)=\frac{1}{\sqrt{2\left[(e_1+e_3-2e_2)(s-e_2)+2(e_1-e_2)(e_2-e_3)
+(e_1-e_3)|s-e_2|\right]}}\\
\times F\!\left(\frac{1}{2},\frac{1}{2},1;
\frac{(e_1+e_3-2e_2)(s-e_2)+2(e_1-e_2)(e_2-e_3)-(e_1-e_3)|s-e_2|}
{(e_1+e_3-2e_2)(s-e_2)+2(e_1-e_2)(e_2-e_3)+(e_1-e_3)|s-e_2|}\right),
\end{multline*}
that is,
\begin{equation*}
\begin{split}
\rho_Q(s)&=\frac{1}{2\pi}\int_{e_2}^{e_1}
\frac{dx}{\sqrt{(e_1-x)(x-e_2)(x-e_3)(x-s)}}\text{ when } e_3<s<e_2,\\
\rho_Q(s)&=\frac{1}{2\pi}\int_{e_3}^{e_2}
\frac{dx}{\sqrt{(e_1-x)(e_2-x)(x-e_3)(s-x)}}
\text{ when } e_2<s<e_1.
\end{split}
\end{equation*}
\end{itemize}
\end{theorem}

\begin{theorem}\label{th:eq}
The density function $\rho_Q$ defined in
Theorem~\ref{th:main} satisfies the following Heun differential
equation
\begin{equation}\label{eq:Heun}
8Q(s)\rho''_Q(s)+8Q'(s)\rho'_Q(s)+Q''(s)\rho_Q(s)=0,
\end{equation}
where $Q(s)=(s-e_1)(s-e_2)(s-e_3)$. Both indices of this equation at the
finite regular singularities $e_1,e_2,e_3$ vanish while its indices at
$\infty$ equal $\frac{1}{2}$ and $\frac{3}{2}$.
\end{theorem}

\begin{corollary}\label{cor:main}
The root-counting measures for each of the eight normalized spectral
polynomials $R_n^{\kappa_1,\kappa_2,\kappa_3}(n(n+1)t)$ as well as that of the
normalized spectral polynomial $R_n(n(n+1)t)$ converge to the measure $\mu_Q$
defined in Theorem~\ref{th:main}, see Figure~\ref{fig1} below.
\end{corollary}

\begin{figure}[!htb]
\centerline{\hbox{\epsfysize=7cm\epsfbox{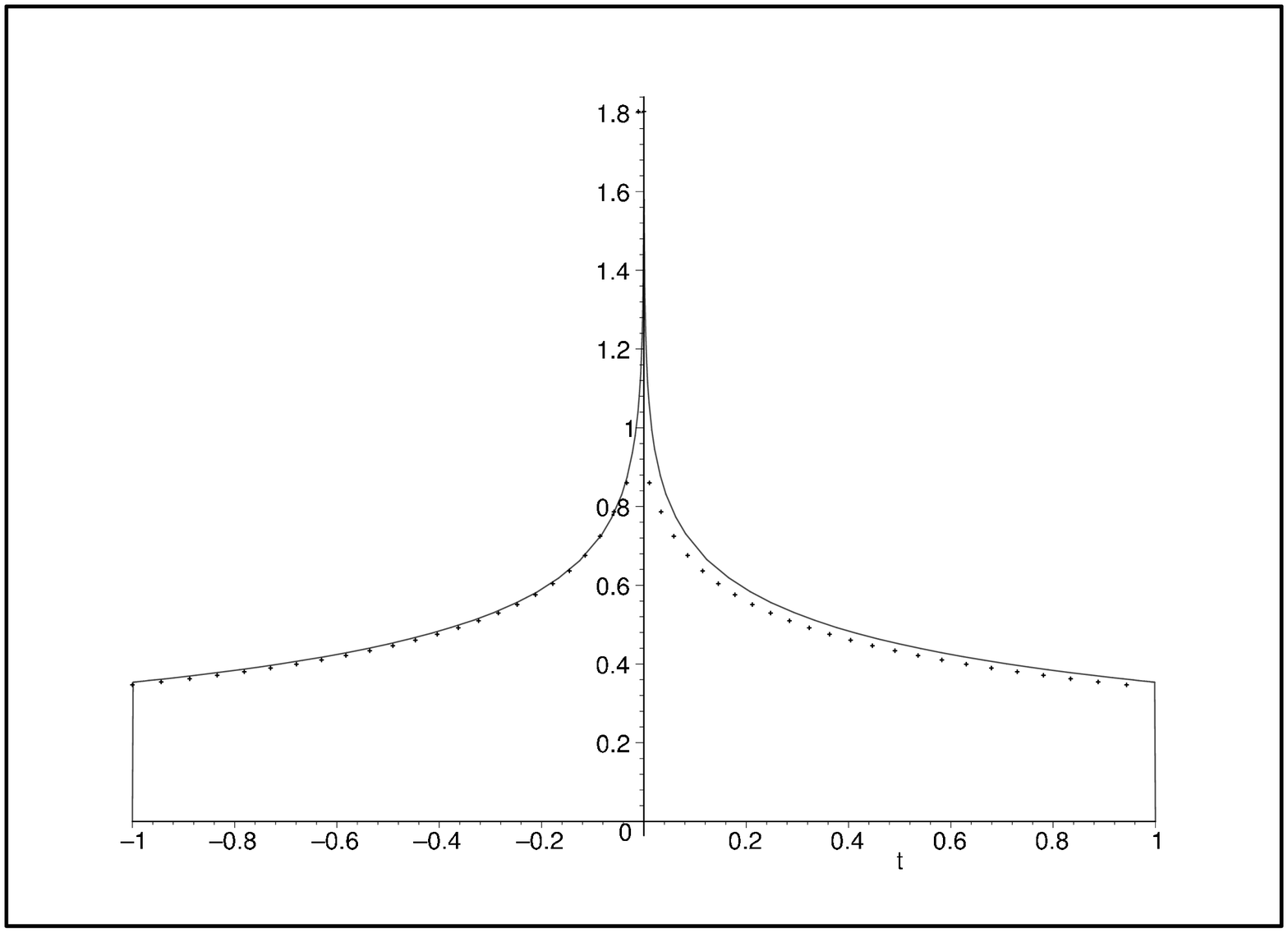}}}
\caption{Comparison of the theoretical density with the numerical density of
the measure $\mu_Q$ for $Q(z)=z^3-z$.}
\label{fig1}
\end{figure}

\begin{remark}\label{r-wei}
Note that equation \eqref{eq:myLame} is usually  lifted to the elliptic
curve $\Psi$ defined by $y^2=Q(z)$ and that on $\Psi$ this equation
takes the so-called Weierstrassian form
$$\left\{\frac{d^2}{du^2}-[n(n+1)\wp(u)+E]\right\}S(u)=0,$$
where $u$ is the canonical coordinate on the universal covering of
$\Psi$ given by
$$u=\int_z^\infty \frac{dw}{\sqrt{Q(w)}}.$$
If $Q(z)$ has all real roots then $n(n+1)\wp(u)$ becomes a finite gap real
periodic potential on $\bR$, see \cite{In1}. Depending on the parity on $n$
four of the eight spectral polynomials $R_n^{\kappa_1,\kappa_2,\kappa_3}(E)$
determine the ends of bands of the spectrum of the corresponding
Schr\"odinger operator on $\bR$. Hopefully our results will find applications
to the spectral theory of finite gap potentials.
\end{remark}

To end this introduction let us mention that in the present context one can
actually ``guess'' the last two formulas in Theorem~\ref{th:main} (iv) by
WKB-type
considerations, as we were kindly informed by K.~Takemura in the final stages
of this work. However, such arguments fail to provide an accurate mathematical
proof. By contrast, our methods use solely rigorous results involving
orthogonal polynomials, elliptic integrals and hypergeometric functions, and
we are currently unaware of any possible ``shortcuts'' in this set-up.
Finally, we should also note that so far neither our methods nor
WKB-type ``guesses'' could provide an answer to these questions in
more general situations, see \S 3 below.

\begin{ack}
We are grateful to M.~Tater from the Institute of Nuclear Physics, Czech
Academy of Sciences,
for interesting discussions on this subject during his visit to Stockholm in
November 2006. We would especially like to thank K.~Takemura from
Yokohama City University for numerous illuminating discussions prior and
during his visit to Stockholm in September 2007. We would also like to
thank the anonymous referee for his insightful comments and suggestions.
\end{ack}

\section{Proofs}

\begin{proof}[Proof of Theorem~\ref{th:main}]
Theorem~\ref{th:main} essentially follows from the main result of \cite{KvA}
after some amount of work.  First we express the polynomial $Sp_n(t)$ as the
characteristic polynomial of a certain matrix. In order to make this matrix
tridiagonal so as to simplify the calculations we assume wlog that
$$Q(z)=(z-e_3)z(z-e_1)=z^3+vz^2+wz$$
with $e_3<e_2=0<e_1$. Set
$$T=(z^3+vz^2+wz)\frac{d^2}{dz^2}+(\al z^2+\be z+\ga)\frac{d}{dz}
-\theta_n(z-t),$$
where $v,w,\al,\be,\ga$ are fixed constants and $\theta_n, t$ are variables.
Assuming  that $S(z)=a_0z^n+a_1z^{n-1}+\ldots+a_n$ with undetermined
     coefficients  $a_i$, $0\le i\le n$, we are looking for the values of
$\theta_n, t$ and $a_i$, $0\le i\le n$, such that $T(S(z))=0$. Note that
     $T(S(z))$ is in general a polynomial of degree $n+1$ whose leading
coefficient equals $a_0[n(n-1)+\al n -\theta_n]$. To get a non-trivial
solution we therefore set
$$\theta_n=n(n-1+\al).$$
Straightforward computations show that the coefficients of
     the successive powers  $z^n, z^{n-1},\ldots, z^0$ in
     $T(S(z))$ can be expressed in the  form of a matrix product
$M_nA$, where $A=(a_0,a_1,\ldots,a_n)^T$ and  $M_n$ is the following
tridiagonal $(n+1)\times(n+1)$ matrix
$$
M_n:=\begin{pmatrix}
t-\xi_{n,1}&\al_{n,2}&0&0&\cdots&0&0\\
\ga_{n,2}&t-\xi_{n,2}&\al_{n,3}&0&\cdots&0&0\\
0&\ga_{n,3}&t-\xi_{n,3}&\al_{n,4}&\cdots&0&0\\
\vdots&\vdots&\ddots&\ddots&\ddots&\vdots&\vdots\\
0&0&0&\ddots&\ddots&\al_{n,n}&0\\
0&0&0&\cdots&\ga_{n,n}&t-\xi_{n,n}&\al_{n,n+1}\\
0&0&0&\cdots&0&\ga_{n,n+1}&t-\xi_{n,n+1}
\end{pmatrix}
$$
with
\begin{equation}\label{eq:extra}
\begin{split}
\xi_{n,i}&=-\frac{v(n-i)(n-i+1)+\be (n-i+1)}{\theta_n},
\quad i\in\{1,\ldots,n+1\}, \\
\al_{n,i}&=\frac{(n-i)(n-i+1)+\al(n-i+1)}{\theta_n}-1,
\quad i\in\{2,\ldots,n+1\}, \\
\ga_{n,i}&=\frac{w(n-i+1)(n-i+2)+\ga (n-i+2)}{\theta_n},
\quad i\in\{2,\ldots,n+1\}.
\end{split}
\end{equation}

A similar matrix can be found in \cite{He} and also  in \cite{Tu}.
The matrix $M_n$ depends linearly on the indeterminate $t$ which appears
only on its main diagonal. If the linear system $M_nA=0$ is to have
a nontrivial solution $A=(a_0,a_1,...,a_n)^T$ the determinant
of $M_n$ has to vanish. This gives the polynomial equation
$$Sp_n(t)=\det(M_n)=0.$$

The sequence of polynomials $\{Sp_n(t)\}_{n\in\bZ_+}$ does
not seem to satisfy any
reasonable recurrence relation. In order to overcome this difficulty and to
be able to use  the  technique of $3$-term recurrence relations  with
variable coefficients (which is applicable since $M_n$ is tridiagonal) we
extend the above polynomial sequence by introducing an additional parameter.
Namely, define
$$Sp_{n,i}(t)=\det M_{n,i},\quad i\in\{1,\ldots,n+1\},$$
where $M_{n,i}$ is the upper $i\times i$ principal submatrix of $M_n$.
One can easily check (see, e.g., \cite[p.~20]{Ar})
that the following $3$-term relation holds
\begin{equation}\label{eq:3term}
Sp_{n,i}(t)=(t-\xi_{n,i}) Sp_{n,i-1}(t) - \psi_{n,i} Sp_{n,i-2}(t),
\quad i\in\{1,\ldots,n+1\},
\end{equation}
where $\xi_{n,i}$ is as in \eqref{eq:extra} and
\begin{equation}\label{eq:coeffs}
\psi_{n,i}=\al_{n,i}\ga_{n,i},\quad i\in\{2,\ldots,n+1\}.
\end{equation}
Here we used the (standard) initial conditions
$Sp_{n,0}(t)=1$, $Sp_{n,-1}(t)=0$.  It is well-known that if all
$\xi_{n,i}$'s are real and all $\psi_{n,i}$'s are  positive then the
polynomials $Sp_{n,i}(t)$, $i\in\{0,\ldots, n+1\}$,
form a (finite) sequence of orthogonal polynomials.  In particular, all their
roots are real. Under the assumptions of Theorem~\ref{th:main}
  the reality of $\xi_{n,i}$ is obvious.
Assuming that the positivity of $\psi_{n,i}$ is also settled
(see Lemma~\ref{lm:posit} below) let us complete the proof of
Theorem~\ref{th:main}.
For this we invoke \cite[Theorem 1.4]{KvA} which translated in our notation
claims that if there exist two continuous functions $\xi(\tau)$ and
$\psi(\tau)$, $\tau\in [0,1]$, such that
$$\lim_{i/(n+1)\to \tau} \xi_{i,n}=\xi(\tau),\quad
\lim_{i/(n+1)\to \tau} \psi_{i,n}=\psi(\tau),\quad
\psi(\tau)\ge 0\,\,\,\,\forall\tau\in [0,1],$$
then the density of the asymptotic root-counting measure of the
polynomial sequence $\{Sp_n(t)\}_{n\in\bZ_+}=\{Sp_{n,n+1}(t)\}_{n\in\bZ_+}$
is given by
$$\int_0^1
\omega_{\left[\xi(\tau)-2\sqrt{\psi(\tau)},
\xi(\tau)+2\sqrt{\psi(\tau)}\right]}d\tau,$$
where for any $x<y$ one has
$$
\omega_{[x,y]}(s)=\begin{cases}
\frac {1}{\pi\sqrt{(y-s)(s-x)}} \text{ if } s\in [x,y],\\
0 \text{ \phantom{xxxxxxxxx}   otherwise.}
\end{cases}
$$

 From the explicit formulas for $\xi_{n,i}$ and $\psi_{n,i}$
(see \eqref{eq:extra} and \eqref{eq:coeffs}) one easily gets
\begin{equation*}
\begin{split}
\xi(\tau)&=\lim_{i/(n+1)\to \tau} \xi_{i,n}=-v(1-\tau)^2, \\
\psi(\tau)&=\lim_{i/(n+1)\to \tau} \psi_{i,n}=-w(1-(1-\tau)^2)(1-\tau)^2.
\end{split}
\end{equation*}
Notice that the above limits are independent of the coefficients
$\al,\be,\ga$ and that by the assumption on $Q(z)$ made at the beginning
of this section one has $w=e_3e_1<0$, which in its turn implies
that $\psi(\tau)\ge 0$ for $\tau\in [0,1]$. The required density
$\rho_Q(s)$ is therefore given by
\begin{equation*}
\begin{split}
\rho_Q(s)&=\int_0^1
\frac{d\tau}{\pi\sqrt{\left( \xi(\tau)+2\sqrt{\psi(\tau)}-s\right)
\left(s-\xi(\tau)+2\sqrt{\psi(\tau)}\right)}^+}\\
&=\int_0^1\frac{d\tau}{\pi\sqrt{-4w\left(1-(1-\tau)^2\right)(1-\tau)^2
-\left(v(1-\tau)^2+s\right)^2}^+},
\end{split}
\end{equation*}
where $\sqrt{\cdot}^+$ is meant to remind that the integrand vanishes
whenever the expression under the square root  becomes negative. Introducing
$\sqrt{\nu}=1-\tau$ we get
\begin{equation*}
\rho_Q(s)=
\frac{1}{2\pi}\int_0^1
\frac{d\nu}{\sqrt{\nu[(4w-v^2)\nu^2-(4w+2vs)\nu-s^2]}^+}.
\end{equation*}
In order to get rid of $\sqrt{\cdot}^+$ we rewrite
\begin{equation}\label{eq:newrho}
\rho_Q(s)=\frac{1}{2\pi\sqrt{v^2-4w}}\int_{\nu_{min}(s)}^{\nu_{max}(s)}
\frac{d\nu}{\sqrt{\nu(\nu-\nu_{min}(s))(\nu_{max}(s)-\nu)}},
\end{equation}
where $\nu_{min}(s)$ and $\nu_{max}(s)$ are the minimal and maximal roots of
of the equation
$$(4w-v^2)\nu^2-(4w+2vs)\nu-s^2=0$$
with respect to $\nu$, that is,
\begin{equation*}
\begin{cases}
\nu_{min}(s)=-(v^2-4w)^{-1}[2w+vs+2\sqrt{w(w+vs+s^2)}],\\
\nu_{max}(s)=-(v^2-4w)^{-1}[2w+vs-2\sqrt{w(w+vs+s^2)}].
\end{cases}
\end{equation*}

A few remarks are in order at this stage. First, by the real-rootedness of
$Q(z)$ one has $v^2-4w=(e_1-e_3)^2>0$. Second, we claim that
$$\left[\xi(\tau)-2\sqrt{\psi(\tau)}, \xi(\tau)+2\sqrt{\psi(\tau)}\right]
\subseteq [e_3,e_1]$$
for $\tau\in [0,1]$. Indeed,
\begin{multline*}
e_1-\xi(\tau)-2\sqrt{\psi(\tau)}=e_1+v(1-\tau)^2
-2\sqrt{-w(1-(1-\tau)^2)(1-\tau)^2}\\
=e_1(1-(1-\tau)^2)-2\sqrt{e_1(1-(1-\tau)^2)}\sqrt{-e_3(1-\tau)^2}
-e_3(1-\tau)^2\\
=\left[\sqrt{e_1(1-(1-\tau)^2)}-\sqrt{-e_3(1-\tau)^2}\right]^2\ge 0
\end{multline*}
and similarly
\begin{multline*}
\xi(\tau)-2\sqrt{\psi(\tau)}-e_3
=-v(1-\tau)^2-2\sqrt{-w(1-(1-\tau)^2)(1-\tau)^2}-e_3\\
=e_1(1-\tau)^2-2\sqrt{e_1(1-\tau)^2}\sqrt{-e_3(1-(1-\tau)^2)}
-e_3(1-(1-\tau)^2)\\
=\left[\sqrt{e_1(1-\tau)^2}-\sqrt{-e_3(1-(1-\tau)^2)}\right]^2\ge 0.
\end{multline*}
It follows that whenever $s$ is such that the integrand in the above formulas
is non-vanishing one has $e_3\le s\le e_1$ hence
\begin{equation*}
\begin{split}
&w(w+vs+s^2)=-e_3e_1(s-e_3)(e_1-s)\ge 0\text{ and}\\
&2w+vs=e_1(e_3-s)+e_3(e_1-s)<0
\end{split}
\end{equation*}
for all $s$ as above. Therefore
\begin{equation*}
\begin{split}
&\nu_{min}(s)+\nu_{max}(s)=-2(2w+vs)(v^2-4w)>0\text{ and}\\
&\nu_{min}(s)\nu_{max}(s)=s^2(v^2-4w)\ge 0
\end{split}
\end{equation*}
from which we conclude that $\nu_{max}(s)>\nu_{min}(s)\ge 0$ for $e_3<s<e_1$.

Recall the definition and properties of the function $\om$ from
\eqref{eq:om}-\eqref{eq:omega} and note that
\begin{equation}\label{eq:om-ex}
\om(e_1,e_3,s)=\frac{(e_1-e_3)|s|}{e_1(s-e_3)-e_3(e_1-s)}
=-\frac{|s|\sqrt{v^2-4w}}{2w+vs}
\end{equation}
so that $\nu_{min}(s)$ and $\nu_{max}(s)$ may actually be rewritten as follows:
$$\nu_{min}(s)=\frac{|s|\!\left(1-\sqrt{1-\om(e_1,e_3,s)^2}\right)}
{(e_1-e_3)\om(e_1,e_3,s)},\quad
\nu_{max}(s)=\frac{|s|\!\left(1+\sqrt{1-\om(e_1,e_3,s)^2}\right)}
{(e_1-e_3)\om(e_1,e_3,s)}.$$
Using these expressions combined with the fact that for any $x<y$ one has
$$\int_x^y\frac{d\nu}{\sqrt{\nu(\nu-x)(y-\nu)}}
=\frac{2}{\sqrt{x}}\mathbb K\!\left(\!\sqrt{\frac{x-y}{x}}\right)$$
(which readily follows from \eqref{eq:elliptic}) we deduce
from \eqref{eq:newrho} that
\begin{equation}\label{eq:first}
\rho_Q(s)=\frac{1}{\pi}
\sqrt{\frac{1+\sqrt{1-\om(e_1,e_3,s)^2}}{(e_1-e_3)|s|\om(e_1,e_3,s)}}
\mathbb K\!\left(i\sqrt{\frac{2\sqrt{1-\om(e_1,e_3,s)^2}}
{1-\sqrt{1-\om(e_1,e_3,s)^2}}}\right).
\end{equation}
Now from \eqref{eq:elliptic}-\eqref{eq:hyper} and the well-known identities
$\Ga(1)=1$, $\Ga(\frac{1}{2})=\sqrt{\pi}$ one easily gets
$$\mathbb K(\ze)=\frac{\pi}{2}F\!\left(\frac{1}{2},\frac{1}{2},1;\ze^2\right)$$
hence
\begin{equation}\label{eq:sec}
\rho_Q(s)=\frac{1}{2}\sqrt{\frac{1+\sqrt{1-\om(e_1,e_3,s)^2}}
{(e_1-e_3)|s|\om(e_1,e_3,s)}}
F\!\left(\frac{1}{2},\frac{1}{2},1;-\frac{2\sqrt{1-\om(e_1,e_3,s)^2}}
{1-\sqrt{1-\om(e_1,e_3,s)^2}}\right)
\end{equation}
by \eqref{eq:first}. It is a remarkable fact due to Kummer and Goursat that
for special choices of the numbers $a,b,c$ the hypergeometric series
$F(a,b,c;z)$ obeys certain quadratic transformation laws. To complete the
proof we need precisely such transformation properties, namely formulas
15.3.19-20 in \cite{AS} with $a=\frac{1}{4}$, $b=\frac{1}{2}$ and $c=1$:
$$
\frac{1}{\sqrt{1-\sqrt{z}}}
F\!\left(\frac{1}{2},\frac{1}{2},1;-\frac{2\sqrt{z}}{1-\sqrt{z}}\right)
=\sqrt{\frac{2}{1+\sqrt{1-z}}}
F\!\left(\frac{1}{2},\frac{1}{2},1;\frac{1-\sqrt{1-z}}{1+\sqrt{1-z}}\right).
$$
The above identity for $z=1-\om(e_1,e_3,s)^2$ together with \eqref{eq:sec}
then yields
\begin{equation}\label{eq:third}
\rho_Q(s)=\sqrt{\frac{\om(e_1,e_3,s)}{2(e_1-e_3)|s|(1+\om(e_1,e_3,s))}}
F\!\left(\frac{1}{2},\frac{1}{2},1;
\frac{1-\om(e_1,e_3,s)}{1+\om(e_1,e_3,s)}\right)
\end{equation}
which by \eqref{eq:om-ex} amounts to
\begin{multline}\label{eq:four}
\rho_Q(s)=\frac{1}{\sqrt{2\left[(e_1+e_3)s-2e_1e_3+(e_1-e_3)|s|\right]}}\\
\times F\!\left(\frac{1}{2},\frac{1}{2},1;
\frac{(e_1+e_3)s-2e_1e_3-(e_1-e_3)|s|}{(e_1+e_3)s-2e_1e_3+(e_1-e_3)|s|}\right).
\end{multline}
To prove the last two formulas of Theorem~\ref{th:main} (iv) we use
\eqref{eq:four} and \eqref{eq:hyper} with $a=b=\frac{1}{2}$ and $c=1$
in order to get an integral representation
of $\rho_Q(s)$ in which one makes the following variable substitutions:
\begin{equation*}
t=-\dfrac{e_3(e_1-x)}{e_1(x-e_3)}\text{ if } e_3<s<0,\,\,
t=-\dfrac{e_1(x-e_3)}{e_3(e_1-x)}\text{ if } 0<s<e_1.
\end{equation*}
The desired expressions are then obtained by straightforward computations.

We thus established all formulas stated in Theorem~\ref{th:main} in the special
case when $e_3<e_2=0<e_1$. The general case reduces to this one by noticing
that if $e_3<e_2<e_1$ and $Q(z)=(z-e_3)(z-e_2)(z-e_1)$ then all the above
arguments may be used
for the polynomial $Q(z+e_2)$. Hence the expressions for $\rho_Q(s)$ given in
Theorem~\ref{th:main} in the general case are obtained from
\eqref{eq:first}-\eqref{eq:four} simply by replacing
$e_1$, $e_3$ and $s$ with $e_1-e_2$, $e_3-e_2$ and $s-e_2$, respectively.
\end{proof}

\begin{lemma}\label{lm:posit}
Let $Q(z)=z^3+vz^2+wz$ and $P(z)=\al z^2+\be z+\ga$ be two polynomials
such that $Q(z)$ has three real distinct roots $e_1>e_2=0>e_3$ and
$$\frac{P(z)}{Q(z)}=\frac{\al_1}{z-e_1}+\frac{\al_2}{z-e_2}
+\frac{\al_3}{z-e_3}$$
with $\al_1, \al_2, \al_3>0$. Then the coefficients $\psi_{n,i}$ defined
in \eqref{eq:coeffs} are all positive.
\end{lemma}

\begin{proof}
One immediately gets $\al>0$, which implies that
$\theta_n>0$ for any positive integer $n$. Recall from \eqref{eq:extra}
and \eqref{eq:coeffs} that $\psi_{n,i}=\al_{n,i}\ga_{n,i}$. Now
$$\al_{n,i}=\theta_n^{-1}\{[(n-i)(n-i+1)-n(n-1)]+\al[(n-i+1)-n]\}$$
and since both $\al$ and $\theta_n$  are positive it follows that
$\al_{n,i}<0$. Let us show that under our assumptions one also has
$\ga_{n,i}<0$. For this it is clearly enough to show that both $w$ and
$\ga$ are negative. Obviously, $w=e_3e_1<0$ and since $\al_1, \al_2,\al_3>0$
the quadratic
polynomial $P(z)$ has two real roots interlacing with $e_3,e_2,e_1$. Therefore
$P(z)$ has one positive and one negative root and positive leading
coefficient hence $\ga=P(0)<0$.
\end{proof}

\begin{proof}[Proof of Theorem~\ref{th:eq}]
In order to deduce the differential equation satisfied by $\rho_Q(s)$ we
note first that the restrictions of $\rho_Q(s)$ to $(e_3,e_2)$ and
$(e_2,e_1)$, respectively, are two branches of the same multi-valued analytic
function (this can be seen e.g.~from the last two expressions in
Theorem~\ref{th:main}). Thus it suffices to derive the linear differential
equation satisfied by $\rho_Q(s)$ restricted to, say, $(e_2,e_1)$.
Specializing formula $(\rm{iv})$ of Theorem~\ref{th:main} to this case we get
$$\rho_Q(s)=\frac{I_Q(s)}{2\pi\sqrt{(e_1-e_2)(s-e_3)}},$$
where
$$I_Q(s)=\int_0^1\frac{dw}{\sqrt{w(1-w)
\left(1-\frac{(e_2-e_3)(e_1-s)}{(e_1-e_2)(s-e_3)}w\right)}}$$
and $s\in (e_2,e_1)$. By \eqref{eq:hyper} we see that up to a constant
factor $I_Q(s)$ is the same as the hypergeometric series
$F\left(\frac{1}{2},\frac{1}{2},1;
\frac{(e_2-e_3)(e_1-s)}{(e_1-e_2)(s-e_3)}\right)$,
which is known to satisfy the following Riemann differential equation
(see, e.g., \cite[\S 15.6]{AS}):
$$I_Q''(s)+\left(\frac{1}{s-e_1}+\frac{1}{s-e_2}\right)\!I_Q'(s)
+\frac{(e_3-e_2)(e_3-e_1)}{4(s-e_3)Q(s)}I_Q(s)=0.$$
Substituting $I_Q(s)=K\sqrt{s-e_3}\rho_Q(s)$, where
$K:=2\pi\sqrt{e_1-e_2}$ is but a constant, in the latter equation we get
after some straightforward algebraic manipulations
\begin{equation}\label{eq:nice}
\rho''_Q(s)+\left(\sum_{i=1}^{3}\frac{1}{s-e_i}\right)\rho'_Q(s)
+\frac{3s-\sum_{i=1}^{3}e_i}{4\prod_{i=1}^{3}(s-e_i)}\rho_Q(s)=0,
\end{equation}
which after multiplication by $8Q(s)$ coincides with the required equation
\eqref{eq:Heun}.
To calculate the indices recall from e.g.~\cite{In2}
that for a  second order linear differential equation
$\rho''(s)+a_1(s)\rho'(s)+a_2(s)\rho(s)=0$ its indicial equation at a
finite regular or regular singular point $\tilde s$ has the form
$$\ze(\ze-1)+\al_1\zeta+\al_2=0,$$
where $\al_1=\lim_{s\to\tilde s}(s-\tilde s)a_1(s)$ and
$\al_2=\lim_{s\to\tilde s}(s-\tilde s)^2a_2(s)$. Thus for \eqref{eq:nice}
the indicial equation at each $e_i$, $i=1,2,3$, has the form $\ze^2=0$,
which implies that both corresponding indices vanish and that a solution to
\eqref{eq:nice} might have a logarithmic singularity at any of these points,
see Remark~\ref{rm:log} below.
The indicial equation of $\rho''(s)+a_1(s)\rho'(s)+a_2(s)\rho(s)=0$
at $\infty$ has the form
$$\zeta(\zeta+1)-\tilde\alpha_1\zeta+\tilde\alpha_2=0,$$
where $\tilde\alpha_1=\lim_{s\to\infty}sa_1(s)$ and
$\tilde\alpha_2=\lim_{s\to\infty}s^2a_2(s)$.
Therefore the indicial equation at $\infty$ for \eqref{eq:nice} is
$\zeta(\zeta+1)-3\zeta+\frac{3}{4}=0$ whose roots are $\frac{1}{2}$ and
$\frac{3}{2}$.
\end{proof}

\begin{remark}\label{rm:log}
It is not difficult to show that
$$\rho_Q(s)\approx \frac{1}{2\pi}\frac{\log\!\left(\frac{16(e_1-e_2)(e_2-e_3)}
{(e_1-e_3)|s-e_2|}\right)}{\sqrt{(e_1-e_2)(e_2-e_3)}}\text{ as }s\to e_2,$$
so that $\rho_Q(s)$ always has a logarithmic singularity at $e_2$.
\end{remark}

\begin{proof}[Proof of Corollary~\ref{cor:main}]
  The main idea of the proof of Corollary~\ref{cor:main} is that the
polynomial part of any Lam\'e solution to \eqref{eq:myLame} itself satisfies
a very similar differential equation. Let us illustrate this in the case of
Lam\'e solutions of type 2. (The other cases can be dealt with in the same
way.) Assume that  $\tilde S(z)$ is a polynomial of
degree $n$ such that $S(z):=(z-e_1)^\frac{1}{2}
(z-e_2)^\frac{1}{2}\tilde S(z)$ solves the equation
$$4Q(z)\!\left[S''(z)+\frac{1}{2}\left(\sum_{i=1}^{3}\frac{1}{z-e_i}\right)\!
S'(z)\right]+V(z)S(z)=0$$
for some linear polynomial $V(z)$. Recall that the (unique) root of each
such $V(z)$ is also a root of the scaled spectral polynomial
$R^{\frac{1}{2},\frac{1}{2},0}_n(n(n+1)t)$. Substituting
$S(z)=(z-e_1)^\frac{1}{2}(z-e_2)^\frac{1}{2}\tilde S(z)$
one gets after some straightforward calculations that $\tilde S(z)$ satisfies
the equation
\begin{equation}\label{eq:modif}
4Q(z)\!\left[\tilde S''(z)+\frac{1}{2}\left(\frac{3}{z-e_1}+\frac{3}{z-e_2}
+\frac{1}{z-e_3}\right)\!\tilde S'(z)\right]+\tilde V(z)\tilde S(z)=0,
\end{equation}
where $\tilde V(z)=V(z)+(z-e_1)-(z-e_2)+4(z-e_3)$. Notice that in the
coefficient in front of $\tilde S'(z)$ the numerators of the first two simple
fractions are increased by $1$ while the remaining one is unchanged.
Thus we are practically in the situation covered by
Theorem~\ref{th:main}. The only (slight) difference is that we are not
considering the asymptotic distribution of the roots of polynomials
$\tilde V(z)$ but that of $V(z)$. However, since
$\tilde V(z)=V(z)+(z-e_1)-(z-e_2)+4(z-e_3)$ and the leading coefficient
$\theta_n$ of $V(z)$ tends to $\infty$ as $n\to \infty$ it follows that
both these families of polynomials actually have the same asymptotic root
distribution. Finally, note that by \eqref{eq:spect}
the scaled spectral polynomial $R_n(n(n+1)t)$ equals the product of four
corresponding polynomials $R^{\kappa_1,\kappa_2,\kappa_3}_n(n(n+1)t)$ and that
all eight polynomial families
$\{R^{\kappa_1,\kappa_2,\kappa_3}_n(n(n+1)t)\}_{n\in\bZ_+}$
have the same asymptotic root distribution when $n\to\infty$.
This implies that the family $\{R_n(n(n+1)t)\}_{n\in\bZ_+}$ itself has
the same limiting root distribution.
\end{proof}

\section{Remarks and Conjectures}

\subsection*{1.} 
\begin{figure}[!htb]
\centerline{\hbox{\epsfysize=4cm\epsfbox{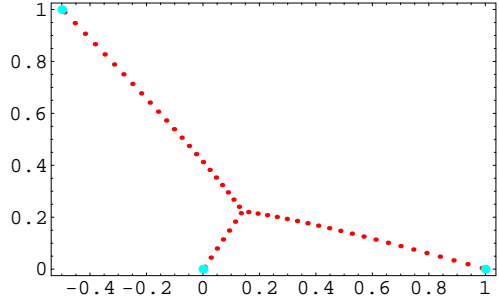}}}
\caption{The roots of the scaled spectral polynomial $Sp_{50}(t)$ for
$Q(z)=z(z-1)\left(z+\frac{1}{2}-i\right)$.}
\label{fig2}
\end{figure}

So far we were unable to extend our method of proving
Theorem~\ref{th:main} to the case of a cubic polynomial $Q(z)$ with complex
roots. The main difficulty comes from the fact that Theorem 1.4 of \cite{KvA}
seems to fail in this case. Indeed, if it were
true then the support of the resulting root-counting measure would be
two-dimensional. However, numerical experiments strongly suggest
that this support is one-dimensional, see Figure~\ref{fig2}.

\begin{problem}
Generalize Theorem~\ref{th:main} to the case of a complex
cubic polynomial $Q(z)$.
\end{problem}

Based on our previous experience (comp.~\cite{B,BBS}) we conjecture that:

\begin{itemize}
\item[(i)] The asymptotic root distribution of the Van Vleck polynomials for
equation \eqref{eq:myGenLame} with complex coefficients is independent of
the $\al_i$'s and depends only on the leading polynomial $Q(z)$;

\item[(ii)] Denoting the above limiting distribution by $\mu_Q$ we claim that
its support is straightened out in the canonical local coordinate
$w(z)=\int_{z_0}^z\frac{dt}{\sqrt{Q(t)}}$.
\end{itemize}

\subsection*{2.} Even more difficulties occur when dealing the more general
Lam\'e equation \eqref{eq:comLame} since in this case Van Vleck polynomials
and their roots can not be found by means of a determinantal equation. They
are in fact related to a more complicated situation when the rank of a
certain non-square matrix is less than the maximal one, see \cite{BBS}.

\begin{problem}
Describe the asymptotic distribution of the roots of all Van Vleck
polynomials for equation \eqref{eq:comLame} when $n\to\infty$.
\end{problem}

An illustration of this asymptotic distribution is given in the next picture.

\begin{figure}[!htb]
\centerline{\hbox{\epsfysize=4cm\epsfbox{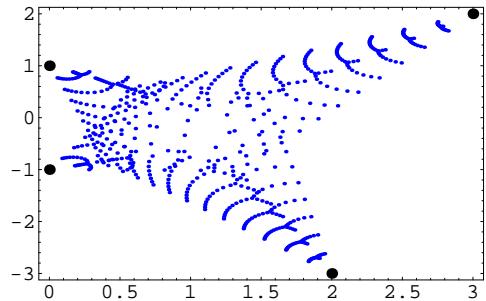}}}
\caption{The union of the roots of 861 quadratic Van Vleck polynomials
corresponding to Stieltjes polynomials of degree $40$ for the classical
Lam\'e equation $Q(z)S''(z)+\frac{Q'(z)}{2}S'(z)+V(z)S(z)=0$ with
$Q(z)=(z^2+1)(z-3i-2)(z+2i-3)$.}
\label{fig3}
\end{figure}

\end{document}